# Family of prime-representing constants: use of the ceiling function

I.A. Weinstein

NANOTECH Centre, Ural Federal University,
Mira street, 19, Ekaterinburg, Russia, 620002

i.a.weinstein@urfu.ru

**Abstract:** The analysis of regularities and randomness in the distribution of prime numbers remains at the research frontiers for many generations of mathematicians from different groups and topical fields. In 2019 D. Fridman et al. [7] have suggested the constant $f_1 = 2.9200509773...$ for generation of the complete sequence of primes with using of a recursive relation for $f_n$ such that the floor $\lfloor f_n \rfloor = p_n$, where $p_n$ is the $n$th prime. Here I present the family of constants $h_n$ ($h_1 = 1.2148208055...$) such that the ceiling function $\lceil h_n \rceil = p_n$. The recursive relation $h_n = \lceil h_{n-1} \rceil (h_{n-1} - \lceil h_{n-1} \rceil + 2)$ generates the sequence of all known prime numbers. I also show constants $h_n$ are irrational.

## 1. Introduction

The analysis of regularities and randomness in the distribution of prime numbers remains at the research frontiers for many generations of mathematicians from different groups and topical fields. The various prime-presenting functions and constants have been proposed and clarified independently [1 – 5, 6 and Refs. in it]. Recently D. Fridman et al. [7] have suggested the constant $f_1 = 2.920050977316...$ for generation of the complete sequence of primes with using of a recursive relation $f_n = \lfloor f_{n-1} \rfloor (f_{n-1} - \lfloor f_{n-1} \rfloor + 1)$ such that the floor of $f_n$ is the $n$th prime, i.e. $\lfloor f_n \rfloor = p_n$. They have also demonstrated $f_1$ is irrational. In this paper I present a family of constants $h_n$ and a recursive relation that generates the complete sequence of primes using the ceiling function approach.

## 2. Recursive relation for $h_n$

**Theorem 1.** *Let $p_n$ denote the nth prime, $n = 1, 2, 3, \ldots$. Then there exists a constant*

$$h_1 = 1.2148208055243337469\ldots$$

*and a sequence*

$$h_n = \lceil h_{n-1} \rceil (h_{n-1} - \lceil h_{n-1} \rceil + 2) \tag{1}$$

*such that the ceiling of $h_n$ is the nth prime, i.e.* $\lceil h_n \rceil = p_n$.

**Proof.** Define $h_1$ as follows:

$$h_1 = \sum_{k=1}^{\infty} \frac{p_k - 2}{\prod_{i=1}^{k-1} p_i}.$$

Considering the oddness of the primes $p_n \geq 3$, well-known Bertrand's postulate [8] can be applied in the following notation $p_n + 1 < p_{n+1} < 2p_n$. Use the formulation for convergence analysis of $h_n$ series:

$$h_1 = (p_1 - 2) + \frac{p_2 - 2}{p_1} + \frac{p_3 - 2}{p_1 p_2} + \frac{p_4 - 2}{p_1 p_2 p_3} + \ldots$$

$$< (p_1 - 2) + \frac{2p_1 - 2}{p_1} + \frac{2p_2 - 2}{p_1 p_2} + \frac{2p_3 - 2}{p_1 p_2 p_3} + \ldots$$

Most terms on the right-hand side of the inequality cancel out, leaving $(p_1 - 2) + 2 = 2$. Since $h_1$ is strictly increasing and bounded, it is convergent: $h_1 = 1.2148208055243337469\ldots$

Define

$$h_n = \sum_{k=n}^{\infty} \frac{p_k - 2}{\prod_{i=n}^{k-1} p_i} = (p_n - 2) + \frac{p_{n+1} - 2}{p_n} + \frac{p_{n+2} - 2}{p_n p_{n+1}} + \frac{p_{n+3} - 2}{p_n p_{n+1} p_{n+2}} + \ldots \tag{2}$$

Obviously, that $h_n = p_{n-1}(h_{n-1} - p_{n-1} + 2)$. Using Bertrand's postulate again, we have

$$h_n > (p_n - 2) + \frac{p_n - 1}{p_n} + \frac{p_{n+1} - 1}{p_n p_{n+1}} + \frac{p_{n+2} - 1}{p_n p_{n+1} p_{n+2}} + \ldots$$

and

$$h_n < (p_n - 2) + \frac{2p_n - 2}{p_n} + \frac{2p_{n+1} - 2}{p_n p_{n+1}} + \frac{2p_{n+2} - 2}{p_n p_{n+1} p_{n+2}} + \ldots$$

Most terms in the inequalities cancel out, leaving $p_n - 1 < h_n < p_n$, and so the ceiling function $\lceil h_n \rceil = p_n$. So, the sequence can now be generated recursively $h_n = \lceil h_{n-1} \rceil (h_{n-1} - \lceil h_{n-1} \rceil + 2)$ ∎

In [7] the corresponding relation for the constants $f_n$ has been proposed

$$f_n = (f_1 - g_{n-1}) \prod_{i=1}^{n-1} p_i = (p_n - 1) + \frac{p_{n+1} - 1}{p_n} + \frac{p_{n+2} - 1}{p_n p_{n+1}} + \frac{p_{n+3} - 1}{p_n p_{n+1} p_{n+2}} + \ldots,$$

where $g_n = \sum_{k=1}^{n} \frac{p_k - 1}{\prod_{i=1}^{k-1} p_i}$ and $f_1 = \lim_{n \to \infty} g_n$.

Obviously that the above expression can be rewritten as similar to Eq. (2):

$$f_n = \sum_{k=n}^{\infty} \frac{p_k - 1}{\prod_{i=n}^{k-1} p_i} \tag{3}$$

Thus, we have two families of prime-representing constants:

$$\lfloor f_n \rfloor = \lceil h_n \rceil = p_n \tag{4}$$

The values of $f_n$ and $h_n$ are given in Table 1 for $n = 1\ldots30$. I've not found constants $h_n$ mentioned anywhere in the independent papers before.

**Table 1.** $f_n$ and $h_n$ families of prime-representing constants

| n | $\lfloor f_n \rfloor = \lceil h_n \rceil = p_n$ | | |
|---|---|---|---|
| | $f_n$ | $h_n$ | $p_n$ |
| 1 | 2.9200509773161347121… | 1.2148208055243337469… | 2 |
| 2 | 3.8401019546322694242… | 2.4296416110486674939… | 3 |
| 3 | 5.5203058638968082726… | 4.2889248331460024817… | 5 |
| 4 | 7.6015293194840413628… | 6.4446241657300124084… | 7 |
| 5 | 11.210705236388289539… | 10.112369160110086858… | 11 |
| 6 | 13.317757600271184934… | 12.236060761210955443… | 13 |
| 7 | 17.130848803525404140… | 16.068789895742420762… | 17 |
| 8 | 19.224429599931870374… | 18.169428227621152959… | 19 |
| 9 | 23.264163538705537112… | 22.219136324801906220… | 23 |
| 10 | 29.075761390227353567… | 28.040135470443843060… | 29 |
| 11 | 31.197080316593253438… | 30.163928642871448728… | 31 |
| 12 | 37.109489814390856584… | 36.081787929014910553… | 37 |
| 13 | 41.051123132461693602… | 40.026153373551690464… | 41 |
| 14 | 43.096048430929437699… | 42.072288315619309040… | 43 |
| 15 | 47.130082529965821066… | 46.108397571630288701… | 47 |
| 16 | 53.113878908393590122… | 52.094685866623568934… | 53 |
| 17 | 59.035582144860276459… | 58.018350931049153499… | 59 |
| 18 | 61.099346546756311072… | 60.082704931900056448… | 61 |
| 19 | 67.060139352134975368… | 66.045000845903443346… | 67 |
| 20 | 71.029336593043349634… | 70.015056675530704168… | 71 |
| 21 | 73.082898106077823987… | 72.069023962679995960… | 73 |
| 22 | 79.051561743681151056… | 78.038749275639705081… | 79 |
| 23 | 83.073377750810933402… | 82.061192775536701425… | 83 |
| 24 | 89.090353317307472383… | 88.079000369546218254… | 89 |
| 25 | 97.041445240365042099… | 96.031032889613424636… | 97 |
| 26 | 101.02018831540908365… | 100.010019029250218968… | 101 |

| 27 | 103.03901985631744826… | 102.02921954272115718… | 103 |
| 28 | 107.01904520069717032… | 106.00961290027918966… | 107 |
| 29 | 109.03783647459722466… | 108.02858032987329405… | 109 |
| 30 | 113.12417573109748793… | 112.11525595618905130… | 113 |

## 3. Irrational constants $h_n$

To prove the irrationality of $h_n$ use the approach from [7].

***Theorem 2.*** The prime-generating constants $h_n$ are irrational for arbitrary $n = 1, 2, 3, \ldots$.

***Proof 1.*** Since $p_n - 1 < h_n < p_n$ for all $n$, we may write $h_n = p_n - r_n$, where $0 < r_n < 1$. Assume $h_n$ is rational such that $h_n = u/w$, where $u, w \in N$. Using the recurrence relation

$$h_{n+1} = p_n(h_n - p_n + 2) = p_n(2 - r_n),$$

it is seen that $wh_n$ is an integer for all $n$. Besides, $r_n \leq (w-1)/w$ for all $n$. Write the above expression as follows:

$$2 - r_n = \frac{h_{n+1}}{p_n} = \frac{p_{n+1}}{p_n} - \frac{r_{n+1}}{p_n}.$$

It is known that $p_{n+1}/p_n \to 1$ when $n \to \infty$. Since the $r_n$ are bounded, the right-hand side tends to 1. In this case $\lim_{n \to \infty} r_n = 1$. It contradicts $r_n \leq (w-1)/w$ for all $n$. ∎

It should be noted, D. Fridman et al. [7] have concluded that $f_1$ is irrational. However, their proof has been conducted for arbitrary $n$. It means that all constants $f_n$ are irrational also.

Further let's demonstrate another way to prove all constants $h_n$ are irrational.

***Proof 2.*** Assume $h_n$ is rational such that $h_n = u/w$, where $u, w \in N$. For arbitrary rational $\frac{a}{b} \neq \frac{u}{w}$ we have $ub - aw \neq 0$. In this case

$$\left| h_n - \frac{a}{b} \right| = \frac{|ub - aw|}{wb} \geq \frac{1}{wb}$$

By Euclid's theorem (there are infinitely many prime numbers) and twin prime conjecture (there are infinitely many twin primes) [8], we can always choose $p_{n+m}$ such that $\frac{p_{n+m}}{2} > w$ and $p_{n+m+1} - p_{n+m} = 2$. Define integer numbers

$$b = \prod_{i=n}^{n+m-1} p_i \quad \text{and} \quad a = \left( p_n - 2 + \frac{p_{n+1} - 2}{p_n} + \frac{p_{n+2} - 2}{p_n p_{n+1}} + \ldots + \frac{p_{n+m+1} - 2}{p_n p_{n+1} \cdots p_{n+m}} \right) \prod_{i=n}^{n+m-1} p_i$$

Then

$$\left| h_n - \frac{a}{b} \right| = \left| \frac{p_{n+m+2} - 2}{p_n p_{n+1} \cdots p_{n+m+1}} + \frac{p_{n+m+3} - 2}{p_n p_{n+1} \cdots p_{n+m+2}} + \ldots \right| =$$

$$= \left| \frac{1}{p_n p_{n+1} \cdots p_{n+m} p_{n+m+1}} \left( p_{n+m+2} - 2 + \frac{p_{n+m+3} - 2}{p_{n+m+2}} + \frac{p_{n+m+4} - 2}{p_{n+m+2} p_{n+m+3}} + \ldots \right) \right| =$$

$$= \left| \frac{1}{p_{n+m} p_{n+m+1} \prod_{i=n}^{n+m-1} p_i} \cdot h_{n+m+2} \right| < \frac{p_{n+m+2}}{p_{n+m} p_{n+m+1} b} < \frac{2}{p_{n+m} b}$$

Thus, it is the contradiction:

$$\left| h_n - \frac{a}{b} \right| \geq \frac{1}{wb} \quad \text{and} \quad \left| h_n - \frac{a}{b} \right| < \frac{1}{wb}$$

It means constants $h_n$ are irrational for arbitrary $n$.  ∎

Application of the floor and ceiling for similar prime-presenting functions is known approach. Mills [1] has shown that there exists a constant $A$ such that $\lfloor A^{3^n} \rfloor$ is prime for all positive integers n. Caldwell and Cheng [9, 10] have calculated the smallest constant $A$ = 1.3063778838…. Elsholtz [11] has given unconditional variants that $\lfloor A^{3^{10n}} \rfloor$ is prime, where A = 1.00536773279814724017…., and $\lfloor A^{3^{13n}} \rfloor$ is prime with A = 3.82499980734391461716…. . Toth [5] has proposed $\lceil B^{3^n} \rceil$ that is also a prime-representing function with constant $B$ = 1.2405547052… [10]. It should be noted the sequence of primes generated by $\lceil B^{c^n} \rceil$ is selective and different than the one for $\lfloor A^{c^n} \rfloor$ with the same value of $c$ and the same starting prime [5]. In our case use of prime-presenting sequences $f_n$ and $h_n$ leads to generating the list of all known primes. Besides any $p_n$ can be taken as a starting point in relations (2) and (3).

## 4. Conclusion

The explicit formula for constants $h_n$, such that the ceiling function $\lceil h_n \rceil = p_n$, has been defined. The corresponding recursive relation has been proposed for generating the complete sequence of known prime numbers. It has been shown that $h_n$ are irrational for all $n$. Two complementary families of prime-representing constants $\lfloor f_n \rfloor = \lceil h_n \rceil = p_n$ have been discussed in the paper by similar manner.